\input amstex
\documentstyle{amsppt}
\magnification=\magstep1 \NoRunningHeads
\topmatter
\title
On the bounded cohomology \\
for ergodic nonsingular  actions \\ of amenable groups
\endtitle

\author
Alexandre I. Danilenko
\endauthor

\email
alexandre.danilenko@gmail.com
\endemail

\address
 Institute for Low Temperature Physics
\& Engineering of National Academy of Sciences of Ukraine, 47 Nauky Ave.,
 Kharkiv, 61103, UKRAINE
\endaddress
\email alexandre.danilenko\@gmail.com
\endemail

\abstract Let $\Gamma$ be an amenable   countable discrete group.
Fix an ergodic free nonsingular  action of $\Gamma$ on a nonatomic standard probability space.
Let $G$ be a 
compactly generated locally compact second countable group such that
the closure of the group of inner automorphisms of $G$ is compact in the natural topology.
It is shown that  there exists a {\it bounded} ergodic $G$-valued cocycle
of $\Gamma$.
\endabstract

\endtopmatter

\document

\head 0. Introduction
\endhead
Let $\Gamma$ be an ergodic infinite countable group of nonsingular transformations of a standard non-atomic probability space $(X,\goth B,\mu)$ and let $G$ be a locally compact second countable group.
Denote by $\Cal M(X,G)$ the set of all measurable  maps from $X$ to $G$.
It is a group under the pointwise multiplication.
A {\it cocycle} $c$ of $\Gamma$ with values in $G$ is a  mapping $c:\Gamma\ni\gamma\mapsto c_\gamma\in \Cal M(X,G)$ satisfying the cocycle equation 
$$
c_{\gamma_1\gamma_2}=c_{\gamma_1}\circ\gamma_2\cdot c_{\gamma_2}
$$
 for all $\gamma_1,\gamma_2\in \Gamma$.
 Denote by $Z^1(\Gamma,G)$ the set of all $G$-valued cocycles of $\Gamma$.
Fix  a left Haar measure $\lambda_G$  on $G$.
Given $c\in Z^1(\Gamma,G)$ and $\gamma\in\Gamma$, we define a nonsingular  transformation $\gamma_c$ of the $\sigma$-finite product space $(X\times G,\mu\times\lambda_G)$ by setting
$\gamma_c(x,g):=(\gamma x, c_\gamma(x)g)$.
Then $(\gamma\beta)_c=\gamma_c\beta_c$ for all $\gamma,\beta\in\Gamma$.
The group $\Gamma_c:=(\gamma_c)_{\gamma\in\Gamma}$ is called {\it the $c$-skew product extension}  of $\Gamma$.
If $\Gamma_c$ is ergodic then $c$ is called {\it ergodic}.
From now on we suppose  that 
 $\Gamma$ (endowed with the discrete topology) and $G$ are both amenable.
 Then it is well known that   there 
 exists an ergodic cocycle of $\Gamma$ with values in $G$.
For the proof of this fact we refer  \cite{Her}, \cite{GoSi1}, \cite{GoSi2}, \cite{Is},  \cite{AaWe1} with application of \cite{Co--We}.

 A subtler question arises naturally: is there  a {\it  bounded} ergodic cocycle $c$ of $\Gamma$ with values in $G$?
 We call  $c$ {\it bounded} if for each $\gamma\in\Gamma$, there is a compact subset $K\subset G$ such that $c_\gamma(x)\in K$ at a.e. $x$. 
 An affirmative answer was obtained recently by Aaronson and Weiss in \cite{AaWe2} in the case where $\Gamma$ is  a $\mu$-preserving freely acting  group isomorphic to $\Bbb Z^n$
 and    $G$ is  a closed subgroup of $\Bbb R^d$.
 (A particular case of this result, when $n=1$, follows from their previous work  \cite{AaWe1}.)
 They also raised a problem of generalization of this result.
 This  is the goal of  the present work.

 Our first observation is that if $\Gamma$ is finitely generated and there is a bounded ergodic $G$-valued cocycle of $\Gamma$ then $G$ is  compactly generated.
 This follows from the fact that  a locally compact second countable group containing a dense compactly generated subgroup  is compactly generated itself \cite{FuSh, Theorem~4}.

Next, we recall that 
\roster
\item"---" $G$ is called an FC-{\it group} if all conjugacy classes in $G$ are finite.
\item"---" $G$ is called an [FC]$^-$-{\it group} if all conjugacy classes in $G$ are relatively compact.
\item"---" $G$ is called an [FD]$^-$-{\it group} if the commutator group in $G$ is relatively compact.
\item"---"  a subset $A$ in $G$ is  {\it normal} if $gAg^{-1}=A$ for each $g\in G$;
\item"---"  
  $G$ is called a  [SIN]-{\it group} if there is a fundamental system of normal neighborhoods of $1_G$.
  \item"---"  
  $G$ is called a [FIA]-{\it group} if $G$ is [FC]$^-$-group and [SIN]-group simultaneously. 
\endroster
Of course, each [FD]$^-$-group is an [FC]$^-$-group.
A locally compact second countable group is an  [FD]$^-$-group if and only if it is an extension of compact group
via an Abelian locally compact group.
A compactly generated [FC]$^-$-group is an [FD]$^-$-group.
It is well known that a locally compact second countable group is a [SIN]-group if and only if 
there is a two-sided invariant metric on $G$ compatible with the topology.
The semidirect product $S:=\Bbb Z\ltimes(\Bbb Z/2\Bbb Z)^\Bbb Z$ with $\Bbb Z$ acting on 
$(\Bbb Z/2\Bbb Z)^\Bbb Z$ by the shift of coordinates is a compactly generated  [FD]$^-$-group which is not
a [SIN]-group.
On the other hand, a countable subgroup $\Bbb Z\ltimes\bigoplus_{n\in\Bbb Z}\Bbb Z/2\Bbb Z$ of $S$ with the discrete topology is a finitely generated [SIN]-group which is not an [FC]$^-$-group.
For a locally compact second countable group $G$, the following claims are equivalent:
\roster 
\item"---" $G$ is a [FIA]-group,
\item"---"  the closure of the group of inner automorphisms of $G$ in the group of continuous automorphisms of $G$ is compact in the natural topology \cite{GrMo},
\item"---" $G$ is a [SIN]-group isomorphic to a direct product of $\Bbb R^n\times L$ for some $n\ge 0$, where a subgroup $L$ contains a compact open normal subgroup $O$ such that the quotient $L/O$ is an FC-group \cite{Wi}.
\endroster
From the last claim we deduce that a [SIN]-group $G$ is a compactly generated [FIA]-group if and only if 
$L/O$ is  a finitely generated FC-group.
The latter happens if and only if the center of $L/O$ has  a subgroup $\Bbb Z^m$, for some $m\ge 0$, which is of finite index in $L/O$.
\comment

\definition{Definition 0.1} We say that $G$ {\it has property $A$} if 
for each neighborhood $U$ of $1_G$ there is a neighborhood $V$ of $1_G$ such that
$gVg^{-1}\subset U$ for each $g\in G$.
\enddefinition

Every Abelian group has property $A$.
Every discrete group has property $A$.
If $K$ is a nontrivial compact metrizable group, $G$ is a countable infinite FC-group  and $\omega$ is  the left Bernoulli action of $G$  on $K^G$ then  the semidirect product $K^G\rtimes_\omega G$ is an [FC]$^{-}$-group without property $A$.

\endcomment
We now state the main result of the work.

\proclaim{Theorem 0.1}
Let $\Gamma$ be an amenable countable discrete group.
Suppose that $\Gamma$ acts nonsingularly on a standard nonatomic probability space and the action is free and ergodic.
If $G$ is a compactly generated  [FIA]-group then there exists a bounded ergodic $G$-valued cocycle of $\Gamma$.
\endproclaim

We deduce Theorem 0.1 from the following more general result.

\proclaim{Theorem 0.2} Let \,$\Gamma$ be as in Theorem~\rom{0.1}.
For  a  [FIA]-group $G$,  fix  a countable symmetric subset $H$ generating a dense subgroup in $G$  and let $H^\bullet:=\bigcup_{g\in G}gHg^{-1}$.
\roster
\item"(i)"
If $\Gamma$ is finitely generated and 
 $\Sigma$ is a finite symmetric generator of $\Gamma$ then
there exists an ergodic cocycle 
$
\alpha=(\alpha_\gamma)_{\gamma\in \Gamma}\in Z^1(\Gamma,G)$ 
such that
$$
\alpha_\sigma(x)\in \{1\}\cup H^\bullet \qquad\text{for a.e. $x\in X$ and $\sigma\in\Sigma.$}
\tag0-1
$$
\item"(ii)"
If $\Gamma$ is not finitely generated and 
 $\Sigma$ is an infinite symmetric generator of $\Gamma$ then
there exist a family $(K_\sigma)_{\sigma\in\Sigma}$ of finite subsets in $G$
and
 an ergodic cocycle 
$
\alpha=(\alpha_\gamma)_{\gamma\in \Gamma}\in Z^1(\Gamma,G)$ such that  
$$
\alpha_\sigma (x)\in K_\sigma\cup H^\bullet \qquad\text{at a.e. $x\in X$ for each $\sigma\in\Sigma$.}
$$
\endroster
\endproclaim

While proving the main result we try to adjust some ideas from \cite{AaWe2} to the more general case under consideration. 
For instance, the sought-for ergodic cocycle appears as the limit of a specially selected sequence of coboundaries.
The pointwise ergodic theorem  was an essential tool in \cite{AaWe2}.
The main obstacle now is that 
the nonsingular counterpart of this theorem (RET,  ratio ergodic theorem) does not hold for nonsingular actions of  arbitrary  amenable groups \cite{Ho}. 
Neither the  weak form of RET from the recent paper \cite{Da} helps in this concrete situation.
To bypass the usage of  RET we utilize  the underlying orbit equivalence relations and their filtrations by finite subrelations.
Secondly, in contrast with the probability preserving case, in our ``nonsingular'' constructions, we have always to take into account  the Radon-Nikodym cocycle (it is trivial in the probability preserving case).

We also record one more corollary from Theorem~0.2 about  existense of {\it norm bounded}  ergodic cocycles 
with values in  groups equipped with a norm.

\proclaim{Corollary 0.3} Let $\Gamma$ be as in Theorem~{\rm 0.1}.
Let $\|.\|$ be a norm on a {[FIA]}-group $G$.
If there is a symmetric countable subset $H\subset G$ such that $\sup_{h\in H}\|h\|<\infty$
and the group generated by $H$ is dense in $G$  then there exists an ergodic cocycle $\alpha\in Z^1(\Gamma,G)$ and a mapping $c:\Gamma\to\Bbb R_+$ such that $\|\alpha(\gamma x,x)\|\le c(\gamma)$ at a.e. $x$ for each $\gamma\in\Gamma$.
\endproclaim

This theorem can be illustrated by an example where $G=\bigoplus_{n>0}\Bbb Z$,
 $\|.\|$ is the $L^\infty$-norm on $G$, i.e. $\|(g_1,g_2,\dots)\|=\max_{n>0}|g_n|$ for each element $(g_1,g_2,\dots)\in G$, and $H$ is the unit ball in $G$.
 We note that $H$ is not compact and $G$ is not compactly generated.
Hence we can not apply Theorem~A.
However, by Corollary~0.3,  norm bounded ergodic cocycles of $\Gamma$ exist.

\head 1. Cocycles of ergodic transformation~groups and equivalence relations
\endhead

\subhead  Cocycles of dynamical systems\endsubhead
Endow $\Cal M(X,G)$ with the topology of convergence in measure.
Then $\Cal M(X,G)$  is a Polish group.
The mapping 
$$
Z^1(X,G)\ni c\mapsto (c_\gamma)_{\gamma\in \Gamma}\in \Cal M(X,G)^\Gamma
$$
 is an embedding of  $Z^1(X,G)$ into the Polish infinite product space $\Cal M(X,G)^\Gamma$.
The image of this embedding is closed.
Hence the  embedding induces a Polish topology on $Z^1(X,G)$.
From now on we denote by dist a metric on $Z^1(X,G)$ compatible with this topology.

Two cocycles $c,d\in Z^1(X,G)$ are called {\it cohomologous} if there is a function $a\in \Cal M(X,G)$ such that $c_\gamma=a\circ\gamma\cdot d_\gamma \cdot a^{-1}$ for all $\Gamma\in\Gamma$.
A cocycle that is cohomologous to the trivial one is called a {\it coboundary}.
The cohomology  relation is an equivalence relation on $Z^1(X,G)$.
If a  cocycle is cohomologous to an ergodic one then it is ergodic.
A coboundary is never ergodic (unless $G$ is a singleton).
For each $f\in\Cal M(X,G)$ and $\gamma\in\Gamma$, we let $\Delta_\gamma f:=f\circ\gamma\cdot f^{-1}$.
Then
 the  map 
 $$
 \Delta f:\Gamma\ni\gamma\mapsto \Delta_\gamma f\in\Cal M(X,G)
 $$ 
 is a coboundary of $\Gamma$.
 Of course, every $\Gamma$-coboundary equals $\Delta f$ for some $f\in\Cal M(X,G)$.
 Moreover, $\Delta a=\Delta b$ for some $a,b\in\Cal M(X,G)$ if and only if  there is $g\in G$ such such $a(x)=b(x) \cdot g$ at a.e. $x$ for some $g\in G$.

It is well known that if $\Gamma$ is amenable (in the discrete topology) then the cohomology class of every cocycle is dense in $Z^1(X,G)$.
In particular, the coboundaries are dense in  $Z^1(X,G)$.
Thus, every cocycle is the  limit of a sequence of coboundaries.

Suppose that $\Gamma$ is generated by a symmetric subset $\Sigma\subset\Gamma$.
It is easy to verify that a sequence $(c^{(n)})_{n=1}^\infty$ of $\Gamma$-cocycles  converges in $Z^1(X,G)$ if and only if the sequence
$(c^{(n)}_\gamma)_{n=1}^\infty$ of functions from $\Cal M(X,G)$ converges in measure for each $\gamma\in\Sigma$.

\subhead Measured equivalence relations and their cocycles
\endsubhead
Denote by $\Cal R$ the $\Gamma$-orbit equivalence relation on $X$.
For $x\in X$, let $\Cal R(x)$ denote the $\Cal R$-class of $x$.
 Then $\Cal R(x)=\{\gamma x\mid \gamma\in\Gamma\}$. 
 A Borel subrelation $\Cal S$ of $\Cal R$ is called {\it conservative} if for each
 subset $A\in\goth B$ of positive measure  and almost every $x\in A$, there is $y\in A$
 such that $y\ne x$ and $(x,y)\in\Cal R$. 
 
For a Borel subrelation $\Cal S$ of $\Cal R$, the {\it full group} $[\Cal S]$ of $\Cal S$ consists  of all $\mu$-nonsingular transformations $\theta$ of $X$ such that  $(\theta x,x)\in\Cal S$  for a.a. $x\in X$.
Of course, $\Gamma\subset[\Cal R]$.
We will use  below the following fact:
\roster
\item"(Fa1)"
If $S$ is conservative then for each subset $A\in\goth B$ of positive measure and $\epsilon>0$,
 there exists
a transformation $\tau\in[\Cal S]$ such that $\tau^2=\text{id}$,  $\tau A=A$,   $\tau x\ne x$
at a.e. $x\in A$, $\tau x=x$ for each $x\not\in A$ and $\Big|\frac{d\mu\circ \tau}{d\mu}(x)-1\Big|<\epsilon$ at a.e. $x\in X$.
\endroster

A Borel subrelation $\Cal S$ of $\Cal R$ is called {\it finite} if $\Cal S(x)<\infty$ for a.e. $x\in X$.
 If $\Gamma$ is amenable then $\Cal R$ is {\it hyperfinite}, i.e. there is an increasing sequence of  finite subrelations $\Cal R_1\subset\Cal R_2\subset\cdots$ of $\Cal R$ such that $\bigcup_{n=1}^\infty\Cal R_n=\Cal R$.
 The sequence $(\Cal R_n)_{n=1}^\infty$ is called a {\it filtration} of $\Cal R$.
There exists a filtration $(\Cal R_n)_{n=1}^\infty$ of $\Cal R$ such that
$\#\Cal R_n(x)=2^n$ at a.e. $x$ for all $n\in\Bbb N$.

A {\it cocycle} of $\Cal R$ with values in $G$ is a Borel map $\alpha:\Cal R\to G$ such that
there is a $\mu$-conull subset $B\subset X$ such that
$$
\alpha(x,y)\alpha(y,z)=\alpha(x,z) \qquad\text{for all $(x,y),(y,z)\in\Cal R\cap(B\times B)$}.
$$
Denote by $Z^1(\Cal R,G)$ the set of all $G$-valued cocycles of $\Cal R$.
Two cocycles $\alpha,\beta\in Z^1(\Cal R,G)$ are called {\it cohomologous} if there exists a function $f\in\Cal M(X,G)$ and  a $\mu$-conull subset $B\subset X$ such that
$\alpha(x,y)=f(x)\beta(x,y)f(y)^{-1}$ for all $(x,y)\in\Cal R\cap (B\times B)$.
A cocycle of $\Cal R$ is a {\it coboundary} if it is cohomologous to the trivial one.
For each $\theta\in [\Cal R]$, we can define a function $\alpha_\theta\in\Cal M(X,G)$ by setting
$\alpha_\theta(x):=\alpha(\theta x,x)$, $x\in X$.
It is well known (and easy to see) that if $\Gamma$ acts freely on $X$ then there is a canonical bijection $\Psi$ of $Z^1(\Cal R,G)$ with $Z^1(\Gamma,G)$.
The bijection is given by the formula $\Psi(\alpha):=(\alpha_\gamma)_{\gamma\in \Gamma}$.
Moreover, two cocycles $\alpha,\beta$ of $\Cal R$ are  cohomologous
if and only if the corresponding cocycles $\Psi(\alpha)$ and $\Psi(\beta)$ of $\Gamma$ are cohomologous.
In particular, $\alpha$ is a coboundary if and only if $\Psi(\alpha)$ is a coboundary.
We transfer dist via $\Psi^{-1}$ to $Z^1(\Cal R,G)$.
Then  $Z^1(\Cal R,G)$ endowed with this metric  is a Polish space.

We single out an important cocycle $\rho_{\Cal R,\mu}\in Z^1(\Cal R,\Bbb R^*_+)$, called {\it the Radon-Nikodym cocycle} of $\Cal R$.
It is well defined by the formula
$$
\rho_{\Cal R,\mu}(x,\gamma x):=\frac{d\mu\circ\gamma}{d\mu}(x), \quad \gamma\in\Gamma,\  x\in X.
$$                    
It is well known (and easy to verify) that $\rho_{\Cal R,\mu}$ is a coboundary if and only if there is a $\sigma$-finite $\Gamma$-invariant measure equivalent to $\mu$.

\subhead Essential values of cocycles\endsubhead
Let $\alpha\in Z^1(\Cal R,G)$.
An element $g\in G$ is called an {\it essential value} of $\alpha$ if 
for each subset $A\in\goth B$ of positive measure and a neighborhood $U$ of $g$ there are
a subset $B\subset A$ of positive measure and 
a transformation $\gamma\in\Gamma$ such that $\gamma B\subset A$ and $\alpha_\gamma(x)\in U$
for all $x\in B$ (\cite{FeMo}, \cite{Sc1}, \cite{Sc2}).
The set $r(\alpha)$ of all essential values of $\alpha$ is a closed subgroup of $G$.
By 
Schmidt's ergodicity criterion for a cocycle, 
$\Psi(\alpha)$ is ergodic if and only if  $r(\alpha)=G$ \cite{Sc2}.

We need a ``nonsingular'' version of  the condition EVC considered by Aaronson and Weiss for probability preserving $\Bbb Z^d$-actions in \cite{AaWe2}.
Let $A\subset X$ be a subset of positive measure, $U$  an open subset in $G$ and $\delta\in(0,1)$.

\definition{Definition 1.1}
 $\alpha$ {\it satisfies }EVC$(A,U, \delta)$
if there are a subset $B\subset A$ and a transformation $\theta\in[\Cal R]$ such that
$\theta B\subset A$, $\mu(B)>\delta\mu(A)$, $\alpha(\theta x,x)\in U$ 
and $\Big|\frac{d\mu\circ\theta}{d\mu}(x)-1\Big|<\delta$ for all $x\in B$.
\enddefinition

The following two facts are routine (cf. \cite{AaWe2}).
Therefore we state them without proof.
\roster
\item"(Fa2)"
The subset of all cocycles of $\Cal R$ satisfying  EVC$(A,U, \delta)$ is open in $Z^1(\Cal R,G)$. \item"(Fa3)"
Fix   a countable base $\Cal U$ of neighborhoods of $1_G$ in $G$, a countable dense subfamily   $\goth B_0$ of $\goth B$ and $g\in G$.
If  for every $U\in\Cal U$ there is a real $\delta=\delta(U,g)>0$ such that 
 $\alpha$ satisfies EVC$(A,Ug, \delta)$ for each $A\in\goth B_0$  then $g\in r(\alpha)$.
\endroster

\head 2. Auxiliary lemma
\endhead

In view of \cite{Co--We},
 we may assume without loss of generality
  that $X=\{0,1\}^\Bbb N$ and $\Cal R$ is {\it the tail equivalence relation} on $X$, i.e. two points $x=(x_j)_{j\in\Bbb N}$ and $y=(y_j)_{j\in\Bbb N}$ 
are equivalent if there is $n>0$ such that $x_j=y_j$ for all $j>n$.
For each $n>0$, we define a  subrelation $\Cal S_n$ of $\Cal R$ by setting $x\sim_{\Cal S_n}y$ if $x_j=y_j$ for all $j>n$.
Then $\Cal S_n$ is a finite subrelation of $\Cal R$, $\#\Cal S_n(x)=2^n$ for each $x\in X$,
$\Cal S_1\subset \Cal S_2\subset\cdots$ and $\bigcup_{n=1}^\infty\Cal S_n=\Cal R$.
Thus, $(\Cal S_n)_{n=1}^\infty$ is a filtration of $\Cal R$.
Replacing, if necessary,  $\mu$ with an equivalent probability measure we may also assume that
\roster
\item"{$(\star)$}"
there is a countable dense subgroup $\Cal L\subset\Bbb R^*$ such that
$\rho_{\Cal R,\mu}(x,y)\in\Cal L$ for all $(x,y)\in\Cal R$.
\endroster
For each $n\in\Bbb N$, we
consider two spaces 
 $$
 X_{n}:=\{0,1\}^{\{1,\dots,n\}}\quad\text{and}\quad X^{n}:=\{0,1\}^{\{n+1,n+2,\dots\}}
 $$
and define two mappings $\pi_n: X\to X_n$ and $\pi^n: X\to X^{n}$ by setting 
 $$
 \pi_n(x):=(x_j)_{j=1}^n,\quad\text{and } \pi^n(x):=(x_j)_{j=n+1}^\infty
 $$
 for each $x:=(x_j)_{j\in\Bbb N}\in X$.
 Then $X$ is isomorphic to the direct product $X_n\times X^n$ via the mapping $\pi_n\times\pi^n$.
 It is easy to see that this mapping  transfers $\Cal R$ to the direct product of the transitive equivalence relation on $X_n$ with the tail equivalence relation $\Cal R^n$ on $X^{n}$.
 Let $\mu^{n}$ stand for the projection of $\mu$ onto $X^{n}$ via  $\pi^n$.
Consider a  disintegration of $\mu$ over $\mu^n$:
$$
\mu=\int_{X^{n}}\mu_{n,y}\otimes\delta_y\,d\mu^{n}(y),
$$
where $(\mu_{n,y})_{y\in X^n}$ is the corresponding system of conditional (probability) measures on $X_n$ and $\delta_y$ is the Kronecker measure on $X^n$ supported at $y$ for each $y\in X^n$.
Since $\mu$ is quasi-invariant  under $\Gamma$, it follows that 
\roster
\item"---"  $\Cal R^n$ is $\mu^n$-nonsingular, i.e. $\mu^n$ is quasiinvariant under each transformation  from the full group $[\Cal R^n]$ and 
\item"---" the measures $\mu_{n,x}$ and $\mu_{n,y}$ are equivalent whenever two points $x,y$ from a $\mu^n$-conull subset of $X^n$ are $\Cal R^n$-equivalent.
\endroster
 Moreover,
 one can verify that
$$
\rho_{\Cal R,\mu}(x,y)=\frac{\mu_{n,y}(\pi_n(y))}{\mu_{n,x}(\pi_n(x))}\cdot \rho_{\Cal R^n,\mu^n}(\pi^n(x),\pi^n(y))
\tag2-1
$$
for all  pairs $(x,y)$ of $\Cal R$-equivalent points belonging to a $\mu$-conull subset of $X$.

Suppose that $G$ is an [FC]$^-$-group.
For  $g\in G$, let $ g^\bullet$ denote the conjugacy class of $g$ in $G$.
Then  the closure of $g^{\bullet}$ is  compact in $G$.
We denote this closure by $K_g$.
Then for each neighborhood $U$ of $1_G$, there is a finite subset $D\subset K_g$ such that
$\bigcup_{d\in D}Ud\supset K_g$.
Since for each $d\in D$, there exists an element  $h(d)\in g^\bullet$ with $d\in Uh(d)$, it follows that
$\bigcup_{d\in D}UUh(d)\supset K_g\supset g^\bullet$.
Hence the positive integer
$$
\Lambda(g,U):=\min\bigg\{\# H\,\bigg|\,\text{$H$ is a finite subset of $g^\bullet$ with $\bigcup_{d\in H}Ud\supset g^\bullet$}\bigg\}
$$ 
is well defined.

\definition{Definition 2.1}
Let $\Cal S$ be a finite subrelation of $\Cal R$, $\Sigma$ a finite subset of $\Gamma$,  $H$ a  subset of $G$ and 
$$
\Cal O(\Cal S,\Sigma):=\{x\in X\mid (\sigma x,x)\not\in\Cal S\text{ for some $\sigma\in\Sigma$}\}.
$$ 
A function $f\in\Cal M(X,G)$ is called
\roster
\item"(i)"
 {\it $\Sigma$-inner on $\Cal S$} if $f(x)=1$
for each $x\in \Cal O(\Cal S,\Sigma)$,
\item"(ii)" 
{\it $(\Sigma,H)$-incremental} if $\Delta_\sigma f(x)\in \{1\}\cup H$ 
at a.e. $x\in X$ for each $\sigma\in\Sigma$.
\endroster
\enddefinition

 Suppose that $\Gamma$ is finitely generated.
 Fix a finite symmetric generating subset   $\Sigma$ of $\Gamma$.
The following lemma is a cornerstone  of   the proof of main results.

\proclaim{Lemma 2.2}
Let $G$ be an $[FC]^-$-group, 
 $H$ a  subset of $G$ and $\epsilon>0$.
 Let   a function $f\in\Cal M (X,G)$ take   finitely many values. 
If  $\Sigma$ is a finite symmetric subset of $\Gamma$ such that  $f$ is  {\it $\Sigma$-inner on $\Cal S_n$} for some $n>0$  and $(\Sigma,H)$-incremental then
for each subset $Z\subset X$ of positive measure,  an element $g\in G$ and  a symmetric normal neighborhood $U$ of $1_G$ in $G$, 
there exist $m>n$, $h\in g^\bullet$,
a  function $\widetilde f\in\Cal M(X,G)$ that takes  finitely many values on $X$, 
a subset $Z_{00}$ of positive measure, a real $\delta=\delta(U,g)>0$ 
and a transformation $\theta\in[\Cal R]$  
such that
\roster
\item"(a)"
$\mu(\Cal O(\Cal S_m,\Sigma))<\epsilon$,
\item"(b)"
$\widetilde f$ is {\it $\Sigma$-inner on $\Cal S_m$} and  $(\Sigma,H\cup\{h,h^{-1}\})$-incremental,   
\item"(c)"
 $Z_{00}\cup \theta Z_{00}\subset Z$, $\mu(Z_{00})>\delta \mu(Z)$,
 $\widetilde f(\theta x)\widetilde f(x)^{-1}\in Ug$ and $\Big|\frac{d\mu\circ\theta}{d\mu}(x)-1\Big|<\epsilon$ for each $x\in Z_{00}$ and 
 \item"(d)"
$
\mu(\{x\in X\mid \Delta_\sigma\widetilde f(x)=\Delta_\sigma f(x)\text{ for each $\sigma\in\Sigma$}\})>1-\epsilon.
$
\endroster
\endproclaim

\demo{Proof} 
First, we find $\epsilon'>0$ such that  $\epsilon'<\epsilon$ and 
$$
\text{
$\sum_{\sigma\in\Sigma}\mu(\sigma M)<\epsilon$ for  every subset $M\in\goth B$ with $\mu(M)<\epsilon'$.}\tag2-2
$$
Next, it
 follows from the definition of $\Lambda(g,U)$ that
there exist a Borel subset $Z_0\subset Z$  and an element $h\in g^\bullet$ such that $\mu(Z_0)\ge
{\mu(Z)}/{\Lambda(g,U)}$ and 
$$
Uh\ni  f(x)^{-1}g f( x)     \quad\text{ for all $x\in Z_0$.}\tag2-3
$$
Without loss of generality we may assume (decreasing $\epsilon$, if necessary)  that
$\epsilon$ is ``considerably'' less   than $\mu(Z_0)$.
Let $G^*:=G\sqcup\{*\}$ stand for the one-point compactification of $G$.
We introduce an auxiliary measurable function $f^*:X\to G^*$ by setting
$$
f^*(x)=\cases
f(x), &\text{if $x\in Z_0$}\\
*,&\text{otherwise}.
\endcases
$$
For each $x\in X$, consider  a  sequence  $(f^*(y),\rho_{\Cal R,\mu}(x,y))_{y\in \Cal S_n(x)}$ of   $2^n$  elements  of  $G^*\times\Cal L$.
We say that two finite sequences of $(G^*\times\Lambda)$-elements are {\it equivalent} if they  coincide up to a permutation.
Denote by $\Cal Z(x)$ the equivalence class of $(f^*(y),\rho_{\Cal R,\mu}(x,y))_{ y\in \Cal S_n(x)}$.
Of course, the mapping $\Cal Z:X\ni x\mapsto\Cal Z(x)$ is $\Cal S_n$-invariant.
Since $f$ takes finitely many values and \thetag{$\star$} holds, the image of $\Cal Z$ is countable.
Therefore
there is a finite partition of $X$ into $\Cal S_n$-invariant measurable subsets  $Y_1,\dots, Y_s, Y_{s+1}$
such that $\min_{1\le j\le s}\mu(Y_j)>0$, $\mu(Y_{s+1})<\epsilon$  and  $\Cal Z$ is constant on $Y_j$ for each $j=1,\dots,s$.
Since $\Cal R_n(x)= X_{n}\times\{\pi^n(x)\}$ for each $x\in X$,
 we obtain a partition of $X^{n}$ into measurable subsets  
$Y_1',\dots, Y_{s+1}'$ such that $Y_j=X_{n}\times Y_j'$ for each $j=1,\dots,s+1$.
It follows that  $\min_{1\le j\le s}\mu^n(Y_j')>0$ and $\mu^n(Y_{s+1}')<\epsilon$.
We note that if two points $a,b$ from $ X^{n}$ belong to the same subset $Y_j'$ for some $j\le s$ then there is a bijection $\theta_{a,b}:X_{n}\to X_{n}$ such that 
$$
\mu_{n,b}\circ\theta_{a,b}=\mu_{n,a}\quad\text{and}\quad
f^*(\theta_{a,b}  t, b)=f^*(t,a)\quad\text{for each }t\in X_{n}.
\tag2-4
$$
Moreover, without loss of generality we may assume that the map 
$$
\bigsqcup_{j=1}^s (Y_j'\times Y_j')\ni (a,b)\mapsto \theta_{a,b}\in \text{Aut\,} X_{n}
$$
is measurable.
By Aut\,$X_n$ we mean the group of bijections  of $X_n$ endowed with the discrete topology.

For $l>n$, let $B_l$ stand for the smallest $\Cal S_n$-invariant subset
containing $\Cal O(\Cal S_l,\Sigma)$.
Then $B_l= X_n\times B_l'$ for some subset $B'_l\subset X^n$.
Since $\Sigma$ is finite and $\Sigma x\subset\Gamma x=\Cal R(x)=\bigcup_{l=1}^\infty\Cal S_l(x)$ for a.e. $x\in X$, it follows
that
$\lim_{l\to\infty}\mu(\Cal O(\Cal S_l,\Sigma))=0$.
This, in turn, implies that $\lim_{l\to\infty}\mu(B_l)=0$ and hence
$\lim_{l\to\infty}\mu^n(B_l')=0$.
For $m>n$, we  consider $X^n$ as the direct product $\{0,1\}^{m-n}\times X^m$.
(We mean the isomorphism $X^{(n)}\ni x=(x_{n+1}, x_{n+2},\dots)\mapsto((x_{n+1},\dots, x_m),(x_{m+1},x_{m+2},\dots)\in \{0,1\}^{m-n}\times X^m$.)
Let  
$$
\mu^n=\int_{X^m}\nu_y\otimes\delta_y\,d\mu^{m}(y)\tag2-5
$$
  denote  the corresponding disintegration of $\mu^n$.
Now we choose 
 $m>n$  large so that 
\roster
\item"(i)"
$\mu^n(B_m')<\epsilon'$ and 
\item"(ii)"
there is a partition $D_1,\dots, D_{s+1}$  of the finite space $\{0,1\}^{m-n}$ such that 
$\mu^n(Y_j'\triangle (D_j\times X^m))<\epsilon\mu^n(Y_j')$ for each $j=1,\dots,s+1$.
\endroster
Let $\Cal P$ denote the set of probability measures on $\{0,1\}^{m-n}$.
Since the mapping  $X^m\ni y\mapsto\nu_y\in\Cal P$ is measurable,
there is  a countable partition $E_1,E_2,\dots$ of $X^m$  such that 
\roster
\item"(iii)"
if $y,y'\in E_j$ for some $j\ge 1$ then $\nu_y\sim\nu_{y'}$ and 
$$
\max\bigg\{\bigg|\frac{\nu_y(z)}{\nu_{y'}(z)}-1\biggl|\,\big|\, z\in\{0,1\}^{m-n}\text{ with }{\nu_y(z)>0}\bigg\}<\epsilon.
$$
\endroster
Since $\Cal R$ is conservative, $\Cal R^m$ is conservative too. 
Therefore, in view of (Fa1), 
 there are 
subsets $E_j'\subset E_j$
$j=1,\dots,s$ 
and a  transformation $\iota\in[\Cal R^m]$ such that $\iota^2=\text{id}$,
$E_j'\cap\iota E_j'=\emptyset$,  
 $E_j=E_j'\sqcup\iota E_j'$ for each $j$ and
$$
\bigg |\frac{d\mu^{m}\circ\iota}{d\mu^{m}}(y)-1\bigg|<\epsilon\quad\text{ at a.e. $y\in 
X^m$}.\tag2-6
$$
We now define a nonsingular transformation $\tau$ of $(X^n,\mu^n)$ by setting $\tau:=\text{id}\times\iota$.
Of course, $\tau\in [\Cal R^n]$ and $\tau^2=\text{id}$.
Let $A:=\bigsqcup_j(\{0,1\}^{m-n}\times E_j')$ and $C:=\tau A$.
Then $A\cap C=\emptyset$ and $A\sqcup C=X^n$.
It follows from  \thetag{2-5}, \thetag{2-6} and (iii) that
$$
\bigg |\frac{d\mu^{n}\circ\tau}{d\mu^{n}}(y)-1\bigg|<3\epsilon\quad\text{ at a.e. $y\in X^n$}.\tag2-7
$$
From this inequality,  (ii)  and the obvious fact that $\tau (D_j\times X^m)=D_j\times X^m$, we deduce that 
$$
\mu^n(Y_j'\triangle\tau Y_j')<4\epsilon\mu^n(Y_j')\qquad\text{for each $j=1,\dots,s+1$}.\tag2-8
$$
We now define  a function $\widetilde f\in\Cal M(X,G)$ by setting
$$
\widetilde f(x)=\cases
1,&\text{if }\pi^n(x)\in B_m'\\
f(x), &\text{if }\pi^n(x)\in A\setminus B_m'\\
f(x)h, &\text{if $\pi^n(x)\in C\setminus B_m'$. }
\endcases
$$
Then $\widetilde f$ takes finitely many values.
We claim that  $\widetilde f$ is $\Sigma$-inner on $\Cal R_m$.
Indeed, if $x\in\Cal O(\Cal S_m,\Sigma)$
then $x\in B_m$ by the definition of $B_m$.
Hence $\pi^{n}(x)\in B_m'$ and $\widetilde f(x)=1$ by the definition of $\widetilde f$, as desired.
Moreover, since $\mu(B_m)<\epsilon$ according to (i),
we obtain (a).

We now verify that $\widetilde f$ is $(\Sigma, H\cup\{h,h^{-1}\})$-incremental.
Take $x\in X$ and $\sigma\in\Sigma$.
Suppose first that $(x,\sigma x)\in\Cal S_n$.
Then 
$\pi_n(x)=\pi_n(\sigma x)$ and
it follows from the definition of $\widetilde f$ that
$$
\Delta_\sigma \widetilde f(x)=
\cases
\Delta_\sigma  f(x)\in H,&\text{if  $x\not\in B_m$}\\
1,&\text{if  $x\in B_m$}.
\endcases
\tag2-9
$$
Consider now the second case:    $(x,\sigma x)\not\in\Cal S_n$.
Then $f(x)=f(\sigma x)=1$ because $f$ is $\Sigma$-inner on $\Cal R_n$ and $\Sigma$ is symmetric. 
Hence 
$$
\Delta_\sigma \widetilde f(x)=
\cases
h, &\text{if   $(\pi^n(\sigma x),\pi^n(x))\in  (C\setminus B_m') \times  (X^n\setminus C)$ }\\
h^{-1}, &\text{if   $(\pi^n(\sigma x),\pi^n(x))\in  (X^n\setminus C)\times(C\setminus B_m') $  }\\
1=\Delta_\sigma  f(x), &\text{otherwise. }\\
\endcases
\tag2-10
$$
Hence $\widetilde f$ is $(\Sigma, H\cup\{h,h^{-1}\})$-incremental.
Thus, (b) is proved.

It follows from~\thetag{2-9} and \thetag{2-10} that  if $\Delta_\sigma\widetilde f(x)\ne\Delta_\sigma f(x)$ then either $x\in B_m$ or 
$(\pi^n(\sigma x),\pi^n(x))\in (C\times (X^n\setminus C))\cup ((X^n\setminus C)\times C)$.
If the latter holds then either  $x\in B_m$ or $\sigma x\in B_m$.
Indeed, if $(\pi^n(\sigma x),\pi^n(x))\in (A\times C)\cup (C\times A)$ then $x\in\Cal O(S_m,\Sigma)\in B_m$
because the subsets $(\pi^n)^{-1}A$ and $(\pi^n)^{-1}C$ are $\Cal S_m$-invariant and disjoint. 
This implies that $\{x\in X\mid \Delta_\sigma\widetilde f(x)\ne\Delta_\sigma f(x)\}\subset B_m\cup \sigma^{-1}B_m$.
Therefore (d) follows from (i) and \thetag{2-2}.

It remains to show (c).
We first let $A':=\bigcup_{j=1}^s(Y_j'\cap\tau Y_j')\setminus(B_m'\cup\tau B_m')$.
It follows from \thetag{2-8}, \thetag{2-7},  (i) and the inequality $\mu^n(Y'_{s+1})<\epsilon$ that
$$
\mu^n(A')>1-7\epsilon.\tag2-11
$$
Now we define a transformation $\theta\in[\Cal R]$ by setting for each $x=(\pi_n(x),\pi^n(x))\in X$,
$$
\theta x:=
\cases
(\theta_{\pi^n(x),\tau\pi^n(x)}\pi_n(x),\tau\pi^n(x)),&\text{if }\pi^n(x)\in A'\\
x, &\text{otherwise. }
\endcases
$$
We see that if $\pi^n(x)\in Y_j'$ for some $x\in X$ and $j\le s$ then $\tau\pi^n( x)\in Y_j'$.
Therefore we deduce from \thetag{2-4}, \thetag{2-1} and \thetag{2-7} that
$$
\bigg|\frac{d\mu\circ\theta}{d\mu}(x)-1\bigg|<3\epsilon \quad\text{and}\quad f^*(\theta x)=f^*(x)
\tag2-12
$$
for a.e. $x\in X$.
The latter equality gives that 
$$
\text{$x\in Z_0$ if and only if $\theta x\in Z_0$ and $f(\theta x)=f(x)$ for each $x\in Z_0$.}\tag2-13
$$
Let $Z_{00}:=\{x\in Z_0\mid \pi^n(x)\in (A\cap  A')\setminus (B_m'\cup\tau B_m')\}$.
Then  
$\theta Z_{00}\cap Z_{00}=\emptyset$.
Since 
$$
Z_0\setminus (Z_{00}\sqcup \theta Z_{00})\subset\{x\in X\mid \pi^n(x)\in (X^n\setminus A')\cup B_m'\cup \tau B_m'\},
$$ 
it follows from \thetag{2-11}, (i)  and the inequality in \thetag{2-12} that $\mu(Z_{00})>\frac12\mu(Z_0)-10\epsilon$.
Next, for each $x\in Z_{00}$, by the definition of $\widetilde f$, \thetag{2-13},  \thetag{2-3} and the normality of $U$,
$$
\widetilde f(\theta x)\widetilde f( x)^{-1}= f(\theta x)h f( x)^{-1}= f(x)h f( x)^{-1}\in f(x)Uf(x)^{-1}g= Ug.
$$
Hence (c) follows if we set  $\delta:=\frac{1}{3\Lambda(g,U)}$.
\qed
\enddemo

\remark{Remark \rom{2.3}}
We note that  if $\Sigma$ generates $\Gamma$ then (d) implies that  dist$(\Delta f,\Delta \widetilde f)$ is ``small'' whenever $\epsilon$ is small enough.
Therefore we can   add  one more  property to the list (a)--(d) in the statement of Lemma~2.2:
\roster
\item"(d)$'$" dist$(\Delta f,\Delta \widetilde f)<\epsilon$.
\endroster
\endremark

\head 3. Proof of the main results
\endhead

In this section we will use the notation from \S2.
In particular, $\Gamma$ is an ergodic nonsingular freely acting amenable countable transformation group
of a standard non-atomic probability space $(X,\goth B,\mu)$ and $(\Cal S_n)_{n=1}^\infty$
is a filtration of the $\Gamma$-orbit equivalence relation.
For $g\in G$ and a normal neighborhood $U$ of $1_G$ in $G$, the real $\delta(U,g)$ is defined in Lemma~2.2.

\demo{Proof of Theorem~\rom{0.2(i)}}
Let $\Sigma$ be a symmetric finite generator of $\Gamma$.
Fix  a countable dense subring  $\goth A$ in $\goth B$ and 
 a countable base  $\goth U$ of normal neighborhoods of $1_G$ in $G$.
 The latter exists because $G$ is a [SIN]-group.
Select a sequence $(A_n,h_n, U_n)_{n=1}^\infty$  of triplets $(A_n,g_n, U_n)\in\goth A\times H\times\goth U$ such that each element of $\goth A\times H\times\goth U$ occurs  in this sequence infinitely many times.
We are going to construct recurrently a sequence $(f_n)_{n=1}^\infty$ of measurable $G$-valued functions on $X$, a decreasing sequence of positive reals $(\epsilon_n)_{n=1}^\infty$
 and an increasing sequence of integers $(k_n)_{n=1}^\infty$ satisfying  the following conditions:
\roster
\item"(i)" $f_n$ takes finitely many values,
\item"(ii)" 
$f_n$ is $\Sigma$-inner on $\Cal S_{k_n}$,
\item"(iii)" 
$f_n$ is  $(\Sigma, H^\bullet)$-incremental,
\item"(iv)" $\mu(\{x\in X\mid \Delta_\sigma f_n(x)=\Delta_\sigma f_{n-1}(x)\text{ for each $\sigma\in\Sigma$}\})>1-\epsilon_n$,
\item"(v)" 
dist$(\Delta f_{n-1},\Delta f_{n})<\epsilon_n$ and 
\item"(vi)"
EVC$(A_n,U_ng_n, \delta(U_n,g_n))$
 holds for $\Delta f_n$.
\endroster
For that we first let $f_1\equiv1_G$.
Suppose now that for  some $n$, we have defined  $k_n$, $\epsilon_n$ and $f_n$ such that
(i)--(vi) hold for $n$.
Utilizing (Fa2) we
choose $\epsilon_{n+1}>0$ such that $\epsilon_{n+1}<\frac12\epsilon_n$ and 
\roster
\item"(vii)"
EVC$(A_n,U_ng_n,\delta(U_n,g_n))$
holds 
for each cocycle $\beta\in Z^1(X,G)$ such that   dist$(\beta,\Delta f_n)<\epsilon_{n+1}$. \endroster
Since (i)--(iii) hold and $G$ is an [FC]$^-$-group, we can apply Lemma~2.2 and Remark~2.3  with $\epsilon_{n+1}$ in place of $\epsilon_n$ to find $k_{n+1}>k_n$ and $f_{n+1}\in\Cal M(X,G)$
such that (i)--(vi) hold with   $n+1$ in place of $n$.
\comment

Let $\Delta f_n$ denote the coboundary generated by $f_n$.
We note that $\mu(C\sqcup A)>1-\epsilon$, $C$ and $A_h$ are $\Cal R$-invariant and $\epsilon$-almost invariant under $\sigma$ for each $h$.
Thus, for most of $x$, either $\{x,\sigma x\}\subset C$ or there is $h\in H$
such that $\{x,\sigma x\}\subset A_h$.

Since for each $\sigma\in \Sigma$ and $n>0$, there is a subset $Y_{\sigma,n}$ of measure greater then $1-\epsilon_n$ such that
$\Delta_\sigma{f_{n+1}}(x)=\Delta_\sigma{f_n}(x)$  for all $x\in Y_{\sigma,n}$.
 It follows from the Borel-Cantelli lemma that there is a map $d:\Sigma\times X\to G$ such that
$d(\sigma,x)=\Delta_\sigma{f_n}(x)$ for a.e. $x$ eventually in $n$.
Thus we obtain a  sequence $(f_n,\epsilon_n,k_n)_{n=1}^\infty$ satisfying (i)--(vi).

\endcomment
Continuing these construction steps infinitely many times, we obtain the entire sequence $(f_n,\epsilon_n,k_n)_{n=1}^\infty$ satisfying (i)--(vi).

Since $\sum_{n=1}^\infty{\epsilon_n}<\infty$, it follows from (v)
that there is  a cocycle $\alpha=(\alpha_\gamma)_{\gamma\in\Gamma}\in Z^1(X,G)$ such that $\Delta f_n\to \alpha$ in dist as $n\to\infty$.
Moreover, (iv) and the Borel-Cantelli lemma yield that 
$\Delta_\sigma f_n(x)=\alpha_\sigma(x)$ for a.e. $x\in X$ eventually in $n$ for each $\sigma\in\Sigma$.
This and (iii) imply~\thetag{0-1}.
Since dist$(\alpha,\Delta f_n)\le \sum_{m> n}\epsilon_m<\epsilon_n$,
it follows from~(vi) and (vii) that EVC$(A_n,U_ng_n,\delta(U_n,g_n))$
holds for $\alpha$ for all $n>0$.
Hence EVC$(A,Ug,\delta(U,g))$
holds for $\alpha$ for all $A\in\goth A$, $g\in H$ and $U\in\goth U$.
It follows from~(Fa3) that $H$ is contained in  $r(\alpha)$.
Since $r(\alpha)$ is a closed subgroup of $G$ and the subgroup generated by $H$ is dense in $G$, we deduce that
 $r(\alpha)=G$, i.e. $\alpha$ is ergodic.
\qed
\enddemo

\demo{Idea  of the proof of Theorem \rom{0.2(ii)}}
Fix an   infinite generator $\Sigma=\{\sigma_j\mid j\in\Bbb N\}$ of $\Gamma$.
For $n>0$, let 
 $\Sigma_n:=\{\sigma_1,\dots,\sigma_n,\sigma_1^{-1},\dots,\sigma_n^{-1}\}$.
Slightly modifying  the proof of Theorem~0.2(i), we construct inductively three sequences $(f_n)_{n=1}^\infty$, $(\epsilon_n)_{n=1}^\infty$ and $(k_n)_{n=1}^\infty$ satisfying (i), (v), (vi) and 
 \roster
 \item"{(ii)$'$}" 
$f_n$ is $\Sigma_n$-inner on $\Cal S_{k_n}$,
\item"(iii)$'$" 
$f_n$ is  $(\{\sigma_j,\sigma_j^{-1}\}, K_j\cup H^\bullet)$-incremental for each $j< n$ 
and $(\{\sigma_n,\sigma_n^{-1}\}, K_n)$-incremental, where $K_j:=F_jF_j^{-1}$ and $F_j\subset G$ is the set of values of $f_j$, $j=1,\dots,n$.
\item"(iv)$'$" $\mu(\{x\in X\mid \Delta_\sigma f_n(x)=\Delta_\sigma f_{n-1}(x)\text{ for each $\sigma\in\Sigma_{n-1}$}\})>1-\epsilon_n$.
 \endroster
 The rest of the proof coincides with the proof of Theorem~0.2(i) almost literally.
 We leave details to the reader.
\qed
\enddemo

\comment

\proclaim{Corollary 2.6}\roster
\item"\rom{(i)}" If \,$\Gamma$ is finitely generated and $G$ is a compactly generated  {[FC]}$^-$-group then there is a bounded ergodic $G$-valued cocycle of $\Gamma$.
\item"\rom{(ii)}" 
If \,$\Gamma$ is not finitely generated and
$G$ is an arbitrary  {[FC]}$^-$-group  then there is a bounded ergodic $G$-valued cocycle of $\Gamma$.
\endroster
\endproclaim

\endcomment

\demo{Proof of  Theorem  \rom{0.1}}
 Since $G$ is compactly generated, there is a symmetric compact subset $K\subset G$  
such that $\bigcup_{n>0}K^n=G$.
Since $G$ is a compactly generated [FC]$^-$-group, it is  an [FD]$^-$-group.
Hence $G$  is an extension of a compact group by an Abelian group.
 It follows that $K^\bullet$ is relatively compact.
Choose a countable symmetric subset $H\subset K$
which is dense in $K$.
Then $H^\bullet$ is relatively compact and the subgroup generated by $H$ is dense in $G$.
It remains to apply Theorem~0.2.
\qed
\enddemo

\enddocument